# Cumulative record times in a Poisson process

CHARLES M. GOLDIE*† and RUDOLF GRÜBEL‡

†Department of Mathematics, Mantell Building, University of Sussex, Brighton BN1 9RF, UK
‡Institut für Mathematische Stochastik, Leibniz Universität Hannover,
Postfach 60 09, D-30060 Hannover, Germany

We obtain a strong law of large numbers and a functional central limit theorem, as $t \to \infty$, for the number of records up to time $t$ and the Lebesgue measure (length) of the subset of the time interval $[0, t]$ during which the Poisson process is in a record lifetime.

## 1. Introduction and results

A Poisson process $N = (N_t)_{t\geq 0}$ with constant intensity $\lambda$ may be regarded as the counting process associated with the partial sums $(S_n)_{n\in\mathbb{N}_0}$ of a sequence $(X_k)_{k\in\mathbb{N}}$ of independent random variables, i.e.

$$N_t = \sup\{n \in \mathbb{N}_0 : S_n \leq t\} \quad \text{with} \;\; S_0 = 0,\; S_n = \sum_{k=1}^{n} X_k \;\; \text{for } n \geq 1.$$

Here each interarrival time $X_k$ has an exponential distribution with parameter $\lambda$,

$$P(X_k \leq x) \;=\; 1 - e^{-\lambda x} \quad \text{for all} \;\; x \geq 0,$$

which we abbreviate to $X_k \sim \mathrm{Exp}(\lambda)$. In renewal-theoretic terms, for which see e.g. [6], [7], $(N_t + 1)_{t\geq 0}$ is the renewal process associated with the lifetimes $(X_k)_{k\in\mathbb{N}}$.

With $1_A$ the indicator function of the set $A$, let $(I_n)_{n\in\mathbb{N}}$ be defined by

$$I_1 = 1, \quad I_n = \prod_{k=1}^{n-1} 1_{\{X_n > X_k\}} \;\; \text{for } n \geq 2,$$

so that $I_n$ indicates whether $X_n$ is a record within the sequence $X_1, X_2, X_3, \ldots$ or not. We are interested in the various aspects of the record structure of $N$, in particular in the two processes $C = (C_t)_{t\geq 0}$ and $W = (W_t)_{t\geq 0}$ given by

$$C_t \;:=\; \#\{1 \leq n \leq N_t : I_n = 1\},$$

*Corresponding author. Email: C.M.Goldie@sussex.ac.uk

and

$$W_t \ := \ \int_0^t I_{N_s+1}\, ds \ = \ \mathrm{Leb}\big(\{0 \le s \le t : \ I_{N_s+1} = 1\}\big)$$

for all $t \ge 0$; we put $C_t = 0$ if $N_t = 0$. In words: $C_t$ is the number of (completed) record inter-arrival times observed up to time $t$, and $W_t$ is the total amount of time up to time $t$ that the Poisson process spends in records. All these variables are defined on some background probability space $(\Omega, \mathcal{A}, P)$. Figure 1 below shows these processes on a finite time interval for a particular $\omega \in \Omega$.

Our first result is a strong law of large numbers for the individual random variables in these processes.

**Theorem 1.1:** *With probability 1, and as* $t \to \infty$*,*

$$\lim_{t\to\infty} \frac{C_t}{\log t} = 1 \quad \textit{and} \quad \lim_{t\to\infty} \frac{W_t}{(\log t)^2} = \frac{1}{2\lambda}\,.$$

The second part can be seen as an instance of the length biasing phenomenon: if we select an index $i$ uniformly at random from $\{1, \dots, n\}$, then the probability that $X_i$ is a record is about $(\log n)/n$ if $n$ is large. However, if we select $t$ uniformly at random from the time interval $[0, n/\lambda]$ (which contains $n$ renewals on average), then the probability that we hit a record renewal is about $(\log n)^2/(2n)$ for $n$ large.

A standard question emerging from Theorem 1.1 is whether we can re-scale the respective differences such that we obtain non-trivial limit distributions as $t \to \infty$. Our method for dealing with this question gives a considerably stronger result: we show that the whole processes converge. To make this precise we use the framework given in the classic [3]. In particular, we regard the processes as random elements of the space $D[0,1]$ of càdlàg functions on the unit interval, endowed with the Skorohod topology and the associated Borel structure; '$\stackrel{d}{\longrightarrow}$' denotes convergence in distribution. We write $Y(t)$ instead of $Y_t$ if this is typographically more appropriate. Then, in order to obtain a converging sequence $(\tilde{Y}_n)_{n\in\mathbb{N}}$ of processes $\tilde{Y}_n = (\tilde{Y}_n(t))_{0\le t\le 1}$ from a given process $Y = (Y(t))_{t\ge 0}$, we put $\tilde{Y}_n(t) = a_n(Y(\phi(n,t)) - f(n,t))$, where, for each $n \in \mathbb{N}$, $a_n \in \mathbb{R}$, $f(n,\cdot) \in D[0,1]$, and $\phi_n : [0,1] \to [0, \tau_n]$ is a time transformation, by which we mean a strictly increasing and continuous function that maps the unit interval onto the initial time segment from 0 to $\tau_n$, where $\tau_n \to \infty$ as $n \to \infty$. From Theorem 1.1 it should be clear that the usual scaling $\phi(n,t) = nt$ would lead to processes with constant paths in some cases; we will instead use $\phi(n,t) = e^{nt} - 1$, so that $\tau_n = e^n - 1$. In particular, we define the re-scaled processes $\tilde{C}_n$ and $\tilde{W}_n$ by

$$\tilde{C}_n(t) \ := \ n^{-1/2}\big(C(e^{nt} - 1) - nt\big)$$

and

$$\tilde{W}_n(t) \ := \ \frac{\lambda}{n^{3/2}}\Big(W(e^{nt} - 1) - \frac{(nt)^2}{2\lambda}\Big),$$

with $0 \le t \le 1$ throughout.

In our functional central limit theorem we consider the two processes simultaneously, which means that we regard $(\tilde{C}_n, \tilde{W}_n)$ as a random element of the product space $D[0,1] \times D[0,1]$, endowed with the product structures (as is well known, this works for the measurability as well as for the topology aspects because of the separability of the Skorohod topology; cf. [3], p. 225).

**Theorem 1.2:** *As* $n \to \infty$*,*

$$(\tilde{C}_n, \tilde{W}_n) \ \stackrel{d}{\longrightarrow} \ Z := (Z^C, Z^W),$$

*where the limit process* $(Z_t^C, Z_t^W)_{0\le t\le 1}$ *can be constructed from a standard Brownian motion* $B = (B_t)_{0\le t\le 1}$

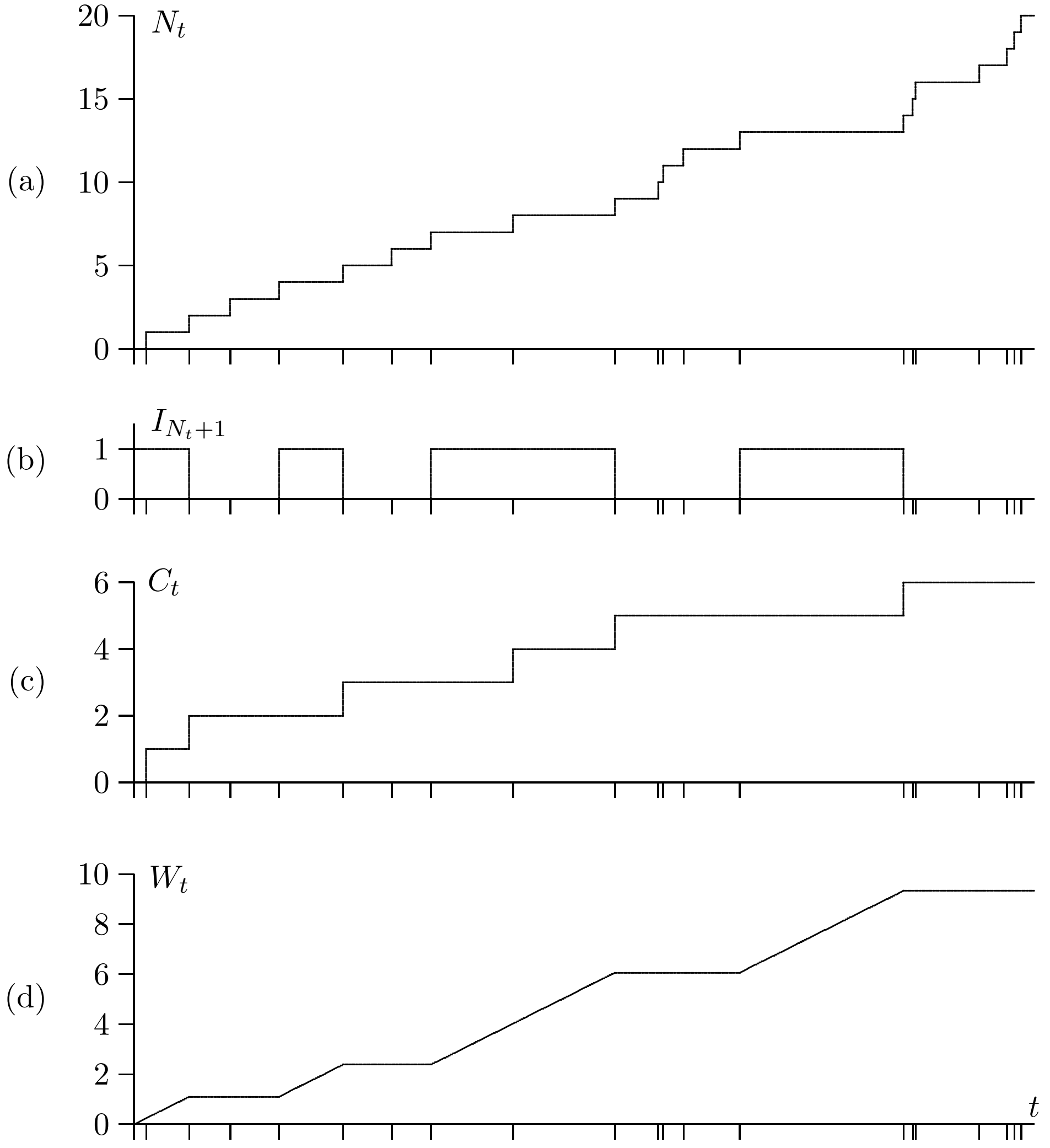

Figure 1. (a) Poisson process, (b) record indicators, (c) number of records, (d) cumulative time in records

*via*

$$Z_t^C = -B_t \quad \textit{and} \quad Z_t^W = \int_0^t B_s \, ds - tB_t, \quad \textit{for } 0 \le t \le 1.$$

Hence, asymptotically and at the level of detail considered here, the second is a simple explicit function of the first.

As will become clear in the proofs, we can replace the discrete parameter $n$ in the definition of the processes by a continuous parameter $\tau$, say, and the convergence holds with $\tau \to \infty$ in the sense that we have convergence for each sequence $(\tau_n)_{n\in\mathbb{N}}$ with $\tau_n \to \infty$ as $n \to \infty$. In particular, we could instead use the scaling $\phi(n,t) = (1+n)^t - 1$, which transforms the unit interval into $[0,n]$.

We consider the $W$-part in more detail. A straightforward calculation (to be given at the end of §2) shows that $Z^W$ is a centred Gaussian process with covariance function $\mathrm{cov}(Z_s^W, Z_t^W) = \min\{s^3, t^3\}/3$, hence $Z^W$ can be seen as a Brownian motion run by a deterministic non-linear clock: the limit process $Z^W$ is equal in distribution to $\big(\tilde{B}(t^3)/\sqrt{3}\big)_{0\le t\le 1}$ where $\tilde{B}(\cdot)$ is another standard Brownian motion.

**Corollary 1.3:** *As $n \to \infty$,*

$$\frac{\lambda\sqrt{3}}{(\log n)^{3/2}}\Big(W\big((n+1)^{t^{1/3}}-1\big)-\frac{t^{2/3}(\log n)^2}{2\lambda}\Big)_{0\le t\le 1} \xrightarrow{d} B,$$

*where $B=(B_t)_{0\le t\le 1}$ denotes standard Brownian motion.*

In particular, with the time transformations $\phi_n(t)=(n+1)^{t^{1/3}}-1$ the initial segments of $W$ have stationary and independent increments in the limit. Of course, Theorem 1.2 also supplies the limit distribution for the individual random variables.

**Corollary 1.4:** *As $t \to \infty$,*

$$\frac{1}{(\log t)^{3/2}}\Big(W_t-\frac{(\log t)^2}{2\lambda}\Big) \xrightarrow{d} W',$$

*where $W'$ has a normal distribution with mean 0 and variance $1/(3\lambda^2)$.*

Apart from asymptotic normality for the individual variables $W_t$ the theorem can also be used to show the asymptotic behaviour of functionals of the process $W$ that depend on the whole history of $W$ up to time $t$, such as in the following corollary.

**Corollary 1.5:** *As $t \to \infty$, and with $B$ as in Corollary 1.3,*

$$\frac{\lambda\sqrt{3}}{(\log t)^{3/2}}\sup_{0\le s\le t}\Big|W_s-\frac{(\log(1+s))^2}{2\lambda}\Big| \xrightarrow{d} \sup_{0\le s\le 1}|B_s|.$$

The limit distribution in Corollary 1.5 is known explicitly; see e.g. [3], p. 80.

Proofs are given in the next section. In §3 we discuss these results and relate them to the literature. We also point out some variants and extensions.

## 2. Proofs

We require some more notation. Let

$$L_1=1,\quad L_{n+1}=\inf\{k>L_n : I_k=1\}\quad\text{for } n\ge 1$$

be the sequence of record times. Then $(R_k)_{k\in\mathbb{N}}$ with $R_i=X_{L_i}$ is the subsequence of record values among the original sequence $(X_k)_{k\in\mathbb{N}}$ of lifetimes. Further, let

$$A_n=\sum_{k=1}^n I_k,\quad T_n=\sum_{k=1}^n R_k$$

be the number of records among the first $n$ lifetimes, and the sum of the first $n$ record values, respectively. Then

$$Z_n=T_{A_n}$$

is the sum of the record values among the first $n$ lifetimes. We will also need to keep in mind that

$$C_t=A_{N_t}.$$

**Proof of Theorem 1.1:** Clearly,

$$Z_{N_t} \le W_t \le Z_{N_t+1} \quad \text{for all } t \ge 0. \tag{1}$$

We need the following two well-known facts from the theory of records of i.i.d. (independent, identically distributed) sequences; see e.g. [1] and [4]. First, it follows from the lack of memory property of exponential distributions that the sequence $(\Delta_k)_{k\in\mathbb{N}}$ of differences of successive record values,

$$\Delta_1 = R_1 \ (= X_1), \quad \Delta_{n+1} = R_{n+1} - R_n \quad \text{for } n \ge 1,$$

is again a sequence of independent random variables with distribution $\mathrm{Exp}(\lambda)$. Secondly, it follows from the permutation invariance of the joint distribution of the $X$-values and the continuity of $\mathrm{Exp}(\lambda)$ that the record indicators $I_n$, for $n \in \mathbb{N}$, are independent and that $P(I_n = 1) = 1/n$ for all $n \in \mathbb{N}$. The second fact implies that $A_n/\log n \to 1$ a.s. (almost surely) as $n \to \infty$; see Lemma 2.1 below. From the first fact and the strong law of large numbers we obtain $R_n/n \to 1/\lambda$ a.s., which together with the elementary fact that

$$\lim_{n\to\infty} a_n = a \quad \Longrightarrow \quad \lim_{n\to\infty} \frac{1}{n^2}\sum_{k=1}^{n} k a_k = \frac{a}{2}$$

for real sequences $(a_n)$ implies that $T_n/n^2 \to 1/(2\lambda)$ a.s. as $n \to \infty$. With probability 1, $A_n \uparrow \infty$, so we can pass to the subsequence and obtain

$$\frac{Z_n}{(\log n)^2} = \frac{T_{A_n}}{A_n^2}\Big(\frac{A_n}{\log n}\Big)^2 \to \frac{1}{2\lambda} \quad \text{a.s. as } n \to \infty. \tag{2}$$

For the second time transformation, from $n$ to $N_t$, we use that $N_t/t \to \lambda$ almost surely as $t \to \infty$; see e.g. p. 107 in [2]. That implies

$$\frac{\log N_t}{\log t} = 1 + \frac{\log(N_t/t)}{\log t} \to 1 \quad \text{a.s. as } t \to \infty. \tag{3}$$

Since $C_t = A_{N_t}$ we gain the first conclusion of the Theorem by replacing $n$ by $N_t$ in the conclusion of Lemma 2.1, and applying (3).

Finally, from (1) we obtain

$$\frac{Z_{N_t}}{(\log N_t)^2}\Big(\frac{\log N_t}{\log t}\Big)^2 \le \frac{W_t}{(\log t)^2} \le \frac{Z_{N_t+1}}{\big(\log(N_t+1)\big)^2}\Big(\frac{\log(N_t+1)}{\log N_t}\Big)^2\Big(\frac{\log N_t}{\log t}\Big)^2.$$

The result for $W_t$ now follows on using (2) and (3), together with the elementary properties of almost-sure convergence. □

**Lemma 2.1:** *$A_n/\log n \to 1$ a.s. as $n \to \infty$.*

**Proof:** Let $B_n = (I_n - n^{-1})/\log(n+1)$, then the $B_n$ are independent, of zero mean, and $E(B_n^2) = (n^{-1} - n^{-2})/(\log(n+1))^2$; consequently $\sum_{n=1}^{\infty} E(B_n^2) < \infty$. It follows, for instance by [14], Theorem 12.2(a), that $\sum_{n=1}^{\infty} B_n$ converges almost surely. By Kronecker's Lemma, $\sum_{k=1}^{n}(I_k - k^{-1})/\log(n+1) \to 0$ a.s. The claimed result follows. □

For the proof of Theorem 1.2 we first observe that by a straightforward re-scaling argument, we may assume $\lambda = 1$. We keep that assumption in force up to the end of Section 2. To start the proof, note that

$$W_t = T_{A_{N_t}} + (t - S_{N_t}) I_{N_t+1}. \tag{4}$$

The distribution of the 'current age' $t - S_{N_t}$ is known explicitly for Poisson processes (it is a truncated exponential) and generally converges as $t \to \infty$ as long as the lifetime distribution has a finite mean. In particular, we have $(t - S_{N_t})I_{N_t+1} = O_P(1)$, and the process version of this bound will allow us to ignore the second term on the right hand side in (4). It is straightforward to obtain asymptotic normality for $(T_n)_{n\in\mathbb{N}}$, but the time change (from $n$ to $A_{N_t}$) requires work. Anscombe's theorem (see e.g. Theorem 7.3.2 in [5]) requires zero-mean summands and then leads to asymptotic normality for $(T_{A_n} - A_n^2/2)/(\log n)^{3/2}$ only. This cannot be combined with a similar statement for $A_n$ unless we have the joint distribution. Hence the proof below provides the joint distribution together with the time change(s). It proceeds by a sequence of lemmas.

We define

$$\lfloor x\rfloor = \max\{n \in \mathbb{Z} : n \le x\},\ [x] = \lfloor x\rfloor + 1,\ \lceil x\rceil = \min\{n \in \mathbb{Z} : n \ge x\}.$$

The first two of these are right-continuous.

**Lemma 2.2:** *The sequence $(R_n)$ of record values satisfies a functional central limit theorem:*

$$\tilde{R}_n := \bigl(n^{-1/2}(R_{[nt]} - nt)\bigr)_{0\le t\le 1} \xrightarrow{d} B = (B_t)_{0\le t\le 1}, \tag{5}$$

*where $B$ denotes Brownian motion.*

**Proof:** We work with the results in Section 17 of [3]. To obtain (5), it is enough to invoke Donsker's theorem together with the result on the record-value differences that was already used in the proof of Theorem 1.1. □

**Lemma 2.3:** *We have*

$$\left(n^{-1/2}\left(n^{-1}T_{[nt]} - \int_0^t R_{[nu]}\,du\right)\right)_{0\le t\le 1} \xrightarrow{P} 0, \tag{6}$$

*and*

$$\bigl(n^{-1/2}(\log L_{[nt]} - R_{[nt]})\bigr)_{0\le t\le 1} \xrightarrow{P} 0. \tag{7}$$

**Proof:** Note that the union of the disjoint intervals $((i-1)/n, i/n]$ for $i = 1, 2, \ldots, [nt]$ is the interval $(0,t]$ plus the interval $(t, [nt]/n]$. Adding the integrals on the intervals $((i-1)/n, i/n]$, and subtracting that on $(t, [nt]/n]$, we obtain

$$\begin{aligned}\int_0^t n^{-1/2}R_{[nu]}\,du &= \frac{1}{n}\sum_{i=1}^{[nt]} n^{-1/2}R_i - \left(\frac{[nt]}{n} - t\right) n^{-1/2}R_{[nt]}\\ &= n^{-3/2}T_{[nt]} - \left(\frac{[nt]}{n} - t\right) n^{-1/2}R_{[nt]}\\ &= n^{-3/2}T_{[nt]} + \theta n^{-3/2}R_{[nt]},\end{aligned}$$

where $\theta = \theta_{n,t}$ lies between $-1$ and 0. Thus

$$\begin{aligned}n^{-1/2}\sup_{0\le t\le 1}\left|n^{-1}T_{[nt]} - \int_0^t R_{[nu]}\,du\right| &\le n^{-3/2}\sup_{0\le t\le 1}|R_{[nt]}|\\ &\le n^{-1}\sup_{0\le t\le 1}|\tilde{R}_n(t)| + n^{-1/2}.\end{aligned} \tag{8}$$

By (5) the supremum converges in distribution to $\sup_{0\le t\le 1}|B_t|$. So the right-hand side of (8) converges to 0 in probability, which proves (6).

For (7) we start from the representation for the record time process given in [4], Theorem 5.8. With $(V_i)_{i\in\mathbb{N}}$ another sequence of independent standard exponentials, also independent of $(R_n)_{n\in\mathbb{N}}$, we have

$$(R_n, L_n)_{n\in\mathbb{N}} \stackrel{d}{=} (R_n, L_n^{(1)})_{n\in\mathbb{N}}$$

where

$$L_n^{(1)} := 1 + \sum_{i=1}^{n-1} \left\lceil \frac{V_i}{-\log(1-e^{-R_i})} \right\rceil,$$

and '$\stackrel{d}{=}$' denotes equality in distribution.

Set $M_1 := \min\{n : R_n \ge \log 2\}$. For $i \ge M_1$ we have $e^{-R_i} \le 1/2$, so the fact that $x \le -\log(1-x) \le 2x$ for $0 \le x \le 1/2$ gives $e^{-R_i} \le -\log(1-e^{-R_i}) \le 2e^{-R_i}$. So for $n > M_1$,

$$\begin{aligned} L_n^{(1)} &= L_{M_1}^{(1)} + \sum_{i=M_1}^{n-1} \left\lceil \frac{V_i}{-\log(1-e^{-R_i})} \right\rceil \\ &= L_{M_1}^{(1)} + \sum_{i=M_1}^{n-1} \lceil V_i \phi_i e^{R_i} \rceil \end{aligned}$$

where $1/2 \le \phi_i \le 1$.

Now let $M_2$ be the least $m \ge M_1$ such that $V_i e^{R_i} \ge 1$ for all $i \ge m$. We need the fact that $V_i e^{R_i}$ tends to infinity a.s., which is left as an exercise for the reader. It implies that $M_2$ is a proper (finite) random variable. For $i \ge M_2$ we have $V_i \phi_i e^{R_i} \ge 1/2$, so the effect of the rounding up given by the $\lceil$ and $\rceil$ brackets can be considered as multiplication by $\psi_i$ where $1 \le \psi_i \le 2$. Thus for $n > M_2$,

$$L_n^{(1)} = L_{M_2}^{(1)} + \sum_{i=M_2}^{n-1} V_i \phi_i \psi_i e^{R_i}.$$

By setting $\phi_i = \psi_i := 1$ for $i < M_2$ and defining $L_n^{(2)} = 1 + \sum_{i=1}^{n-1} V_i \phi_i \psi_i e^{R_i}$ for all $n \ge 1$, we obtain

$$L_n^{(1)} = A + L_n^{(2)} \qquad \text{for } n > M_2,$$

where $A := L_{M_2}^{(1)} - 1 - \sum_{i=1}^{M_2-1} V_i e^{R_i}$.

Next,

$$\log L_n^{(1)} = \log L_n^{(2)} + \log\left(1 + \frac{A}{L_n^{(2)}}\right) \qquad \text{for } n > M_2.$$

For the last term here we use the fact that $L_n^{(2)} > L_{M_2}^{(2)}$ for $n > M_2$, which implies that $|\log(1 + A/L_n^{(2)})|$ increases when we replace $L_n^{(2)}$ in it by $L_{M_2}^{(2)}$ (consider the cases $A > 0$ and $A < 0$ separately). For the other term on the right-hand side we use the fact that $1/2 \le \phi_i\psi_i \le 2$ to deduce that $\log L_n^{(2)}$ differs from $\log\bigl(1 + \sum_{i=1}^{n-1} V_i e^{R_i}\bigr)$ by at most $\pm\log 2$. Therefore

$$(\log L_n - R_n)_{n\in\mathbb{N}} \stackrel{d}{=} \left(\epsilon_n + \log\Bigl(e^{-R_n} + \sum_{i=1}^{n-1} V_i e^{-(R_n - R_i)}\Bigr)\right)_{n\in\mathbb{N}}$$

where

$$|\epsilon_n| \le \log 2 + \left|\log\left(1 + \frac{A}{L^{(2)}_{M_2}}\right)\right| \qquad \text{for } n > M_2.$$

It is clear that $n^{-1/2}\sup_{n\in\mathbb{N}}|\epsilon_n| \xrightarrow{a.s.} 0$ as $n\to\infty$. Formula (7) follows by virtue of Lemma 2.4 below. □

**Lemma 2.4:** *With the $V_i$ as in the previous proof,*

$$n^{-1/2}\max_{k=1,\dots,n}\left|\log\Big(e^{-R_k} + \sum_{i=1}^{k-1} V_i e^{-(R_k-R_i)}\Big)\right| \xrightarrow{P} 0$$

*as $n\to\infty$.*

**Proof:** Choose $\epsilon > 0$, then

$$\begin{aligned} P\left(\log\Big(e^{-R_k} + \sum_{i=1}^{k-1} V_i e^{-(R_k-R_i)}\Big) > \epsilon\sqrt{k}\right) &\le P\Big(1 + \sum_{i=1}^{k-1} V_i > e^{\epsilon\sqrt{k}}\Big) \\ &\le \frac{E(1+\sum_{i=1}^{k-1} V_i)}{e^{\epsilon\sqrt{k}}} = \frac{k}{e^{\epsilon\sqrt{k}}}. \end{aligned}$$

As the latter is summable it follows that, with probability 1, there are at most finitely many $k$ such that $\log(e^{-R_k} + \sum_{i=1}^{k-1} V_i e^{-(R_k-R_i)})$ exceeds $\epsilon\sqrt{k}$, and it follows that

$$n^{-1/2}\max_{k=1,\dots,n}\log\Big(e^{-R_k} + \sum_{i=1}^{k-1} V_i e^{-(R_k-R_i)}\Big) \xrightarrow{P} 0. \tag{9}$$

On the other hand,

$$\begin{aligned} P\left(\log\Big(e^{-R_k} + \sum_{i=1}^{k-1} V_i e^{-(R_k-R_i)}\Big) < -\epsilon\sqrt{k}\right) &= P\Big(e^{-R_k} + \sum_{i=1}^{k-1} V_i e^{-(R_k-R_i)} < e^{-\epsilon\sqrt{k}}\Big) \\ &< P(V_{k-1}e^{-\Delta_k} < e^{-\epsilon\sqrt{k}}) \\ &= P(V_{k-1} < e^{-\epsilon\sqrt{k}}) + P(\Delta_k > \epsilon\sqrt{k} + \log V_{k-1} > 0) \\ &= 1 - \exp(-e^{-\epsilon\sqrt{k}}) + \int_{e^{-\epsilon\sqrt{k}}}^{\infty} e^{-(\epsilon\sqrt{k}+\log x)} e^{-x}\,dx. \end{aligned} \tag{10}$$

In the last right-hand side the term $1-\exp(-e^{-\epsilon\sqrt{k}})$ is $(1+o(1))e^{-\epsilon\sqrt{k}}$ as $k\to\infty$. And the integral is the sum of

$$e^{-\epsilon\sqrt{k}}\int_1^\infty e^{-x}\,\frac{dx}{x} = e^{-\epsilon\sqrt{k}}E_1(1)$$

($E_1(\cdot)$ is the exponential integral), and

$$e^{-\epsilon\sqrt{k}}\int_{e^{-\epsilon\sqrt{k}}}^1 e^{-x}\,\frac{dx}{x} < e^{-\epsilon\sqrt{k}}\int_{e^{-\epsilon\sqrt{k}}}^1 \frac{dx}{x} = \epsilon\sqrt{k}\,e^{-\epsilon\sqrt{k}}.$$

Thus the left-hand side of (10) is bounded by a constant multiple of $\sqrt{k}\,e^{-\epsilon\sqrt{k}}$ and so is summable. Arguing as for (9), it follows that

$$n^{-1/2} \max_{k=1,\ldots,n} -\log\Big(e^{-R_k} + \sum_{i=1}^{k-1} V_i e^{-(R_k - R_i)}\Big) \xrightarrow{P} 0,$$

and with (9) this gives the result. □

Define

$$\tilde{T}_n := \big(n^{-3/2}(T_{[nt]} - (nt)^2/2)\big)_{0 \le t \le 1},$$

the (standardised) record sum process, and

$$\tilde{L}_n := \big(n^{-1/2}(\log L_{[nt]} - nt)\big)_{0 \le t \le 1},$$

the (standardised) record times process.

**Lemma 2.5:** *As $n \to \infty$,*

$$(\tilde{T}_n, \tilde{L}_n) \xrightarrow{d} \Psi(B),$$

*where we define $\Psi : D[0,1] \to D[0,1] \times D[0,1]$ by $\Psi(f) := (\Psi_1(f), f)$, with $\Psi_1 : D[0,1] \to D[0,1]$ in turn defined by $\Psi_1(f)(t) = \int_0^t f(s)\, ds$.*

**Proof:** Clearly, $\Psi_1$ is continuous (and linear), so that $\Psi$ is continuous too. By the continuous mapping theorem, applied to the conclusion of Lemma 2.2, it follows that $\Psi(\tilde{R}_n) \xrightarrow{d} \Psi(B)$. But Lemma 2.3 says that $(\tilde{T}_n, \tilde{L}_n) = \Psi(\tilde{R}_n) + o_P(1)$, so our claimed result follows. □

Now we switch to the 'stochastic clock(s)' $\Phi_n$ given by

$$\Phi_n(t) \;:=\; \frac{1}{n} A_{N(e^{nt}-1)}, \text{ for } 0 \le t \le 1,$$

and deduce joint convergence of the time-changed processes $\tilde{T}_n \circ \Phi_n$ and $\tilde{L}_n \circ \Phi_n$.

**Lemma 2.6:** *As $n \to \infty$,*

$$(\tilde{T}_n \circ \Phi_n, \tilde{L}_n \circ \Phi_n) \xrightarrow{d} \Psi(B).$$

**Proof:** From Lemma 2.1 it follows that

$$\max_{1 \le k \le e^n} \big|A_k - \log(1+k)\big| \;=\; o_P(n). \tag{11}$$

Choose $\epsilon > 0$. Then

$$
\begin{aligned}
&P\left(\sup_{0\le t\le 1} |A_{N(e^{nt}-1)} - \log(1+N(e^{nt}-1))| > \epsilon n\right)\\
&\quad \le P(N(e^n-1) > e^{n+1}) + P\left(N(e^n-1)\le e^{n+1},\ \sup_{0\le t\le 1} |A_{N(e^{nt}-1)} - \log(1+N(e^{nt}-1))| > \epsilon n\right)\\
&\quad \le P(N(e^n-1) > e^{n+1}) + P\left(N(e^n-1)\le e^{n+1},\ \max_{1\le k\le e^{n+1}} |A_k - \log(1+k)| > \epsilon n\right)\\
&\quad \le P(N(e^n-1) > e^{n+1}) + P\left(\max_{1\le k\le e^{n+1}} |A_k - \log(1+k)| > \frac{\epsilon}{2}(n+1)\right),
\end{aligned}
$$

since $\epsilon n \ge \epsilon(n+1)/2$. On the right-hand side, (11) makes the second term tend to 0, while the first term tends to 0 because $N_t/t \to 1$ a.s. We have shown that

$$\sup_{0\le t\le 1} \left|\Phi_n(t) - n^{-1}\log(1+N(e^{nt}-1))\right| = o_P(1),$$

and it follows immediately that

$$\sup_{0\le t\le 1} \left|\Phi_n(t) - t\right| = o_P(1). \tag{12}$$

Using this, we may proceed as in the proof of Theorem 17.3 in [3] to conclude the present proof. □

Inserting $A_{N_t}$ for $n$ in $L_n$ is like time inversion, so we expect that $\log L_{n\Phi_n(t)}$ differs from $nt$ by a negligible amount. That is indeed the conclusion of Lemma 2.12 below, and the intervening results are steps on the road to that end. We allow ourselves to write $A_t$ for non-integer $t$, to be interpreted as $A_{\lfloor t\rfloor}$, and similarly for other random sequences.

**Lemma 2.7:** *As $n \to \infty$,*

$$\sup_{0\le t\le 1} n^{-1/2}\big(\log L(A_{e^{nt}}+1) - \log L(A_{e^{nt}})\big) = o_P(1).$$

**Proof:** For $m \ge k$,

$$P(L_{n+1} > m | L_n = k) = \frac{k}{k+1}\frac{k+1}{k+2}\cdots\frac{m-1}{m} = \frac{k}{m},$$

so

$$P(L_{n+1} > xL_n | L_n = k) = \frac{k}{\lfloor xk\rfloor} \le \frac{1}{x-1},$$

which we use in the form

$$P\left(\frac{L_{n+1}}{L_n} > x\right) \le \frac{1}{x-1} \qquad (x>1).$$

Choose $\epsilon > 0$, then

$$\begin{aligned} P\left(\sup_{0\le t\le 1} \frac{L(A_{e^{nt}}+1)}{L(A_{e^{nt}})} > e^{\epsilon\sqrt{n}}\right) &\le P(A_{e^n} > 2n) + P\left(\max_{k=1,\ldots,2n} \frac{L_{k+1}}{L_k} > e^{\epsilon\sqrt{n}}\right) \\ &\le P(A_{e^n} > 2n) + \sum_{k=1}^{2n} P\left(\frac{L_{k+1}}{L_k} > e^{\epsilon\sqrt{n}}\right) \\ &\le P(A_{e^n} > 2n) + \frac{2n}{e^{\epsilon\sqrt{n}}-1}, \end{aligned}$$

and this tends to 0 as $n \to \infty$ because $A_{e^n}/n \to 1$ almost surely. This establishes Lemma 2.7. □

**Lemma 2.8:** *As* $n \to \infty$,

$$\left(n^{-1/2}\big(\log L(A_{e^{nt}}) - nt\big)\right)_{0\le t\le 1} = o_P(1).$$

**Proof:** For $L_k \le i < L_{k+1}$, $A_i = k$, so $L_{A_i} = L_k$. Therefore $L_{A_t} \le t < L_{A_t+1}$ for all $t$. We use this in the form

$$0 \le \log t - \log L_{A_t} < \log L_{A_t+1} - \log L_{A_t}.$$

Thus

$$\sup_{0\le t\le 1} n^{-1/2}|\log L(A_{e^{nt}}) - nt| \le \sup_{0\le t\le 1} n^{-1/2}\big(\log L(A_{e^{nt}}+1) - \log L(A_{e^{nt}})\big),$$

and the result follows by Lemma 2.7. □

**Lemma 2.9:** *As* $n \to \infty$,

$$\left(n^{-1/2}\big(\log L(A_{e^{nt}-1}+1) - nt\big)\right)_{0\le t\le 1} = o_P(1).$$

**Proof:** We write

$$\sup_{0\le t\le 1}\left|n^{-1/2}\big(\log L(A_{e^{nt}-1}+1) - nt\big)\right| \le \left(\sup_{0\le t<n^{-1}\log 2} + \sup_{n^{-1}\log 2\le t\le 1}\right)\left|n^{-1/2}\big(\log L(A_{e^{nt}-1}+1) - nt\big)\right|.$$

The first supremum on the right equals

$$\sup_{0\le t<n^{-1}\log 2} n^{-1/2}|\log L_1 - nt| = n^{-1/2}\log 2 = o(1).$$

The second is bounded by

$$\begin{aligned} &\sup_{n^{-1}\log 2\le t\le 1} n^{-1/2}\big|\log L(A_{e^{nt}-1}+1) - \log(e^{nt}-1)\big| + \sup_{n^{-1}\log 2\le t\le 1} n^{-1/2}|\log(e^{nt}-1) - nt| \\ &\le \sup_{0\le t\le 1} n^{-1/2}\big|\log L(A_{e^{nt}}+1) - nt\big| + \sup_{n^{-1}\log 2\le t\le 1} n^{-1/2}|\log(1-e^{-nt})| \\ &\le \sup_{0\le t\le 1} n^{-1/2}\big|\log L(A_{e^{nt}}+1) - \log L(A_{e^{nt}})\big| + \sup_{0\le t\le 1} n^{-1/2}\big|\log L(A_{e^{nt}}) - nt\big| + o(1). \end{aligned}$$

The two terms on the final right-hand side are $o_P(1)$ by Lemmas 2.7 and 2.8 respectively. □

**Lemma 2.10:** *As $n \to \infty$,*

$$\sup_{0 \le t \le 1} n^{-1/2} |\log(N(e^{nt}-1)+1) - nt| = o(1) \ \textit{a.s.}. \tag{13}$$

**Proof:** From $N_t \sim t$ a.s., i.e. $N_t/t \to 1$, it follows that $N(e^t-1)+1 \sim e^t$ a.s. Therefore

$$M := \sup_{0 \le t < \infty} |\log(N(e^t-1)+1) - t|$$

is a proper (finite) random variable. The left-hand side of (13) is bounded by $n^{-1/2}M$, so converges to 0 a.s. as $n \to \infty$. □

**Lemma 2.11:** *As $n \to \infty$,*

$$\Big(n^{-1/2}\big(\log L(A_{N(e^{nt}-1)}+1) - nt\big)\Big)_{0 \le t \le 1} = o_P(1).$$

**Proof:** First, because of Lemma 2.10 it suffices to show

$$\Big(n^{-1/2}\big(\log L(A_{N(e^{nt}-1)}+1) - \log(N(e^{nt}-1)+1)\big)\Big)_{0 \le t \le 1} = o_P(1).$$

Choose $\epsilon > 0$. We have

$$\begin{aligned} &P\left(\sup_{0 \le t \le 1}\Big|n^{-1/2}\big(\log L(A_{N(e^{nt}-1)}+1) - \log(N(e^{nt}-1)+1)\big)\Big| > \epsilon\right) \\ &\quad \le P(N(e^n-1)+1 > e^{2n}) \\ &\qquad + P\left(N(e^n-1)+1 \le e^{2n}, \sup_{0 \le t \le 1}\Big|n^{-1/2}\big(\log L(A_{N(e^{nt}-1)}+1) - \log(N(e^{nt}-1)+1)\big)\Big| > \epsilon\right). \end{aligned}$$

Set $s(t) := n^{-1}\log(N(e^{nt}-1)+1)$, then as $t$ varies over its range $0 \le t \le 1$, $s$ varies over the range $0 \le s \le n^{-1}\log(N(e^n-1)+1)$, which is contained in $0 \le s \le 2$ because we are within the event $N(e^n-1)+1 \le e^{2n}$. Thus the above is bounded by

$$\begin{aligned} &P(N(e^n-1)+1 > e^{2n}) + P\left(\sup_{0 \le s \le 2} n^{-1/2}|\log L(A_{e^{ns}-1}+1) - ns| > \epsilon\right) \\ &\quad = P(N(e^n-1)+1 > e^{2n}) + P\left(\sup_{0 \le t \le 1} 2^{1/2} m^{-1/2}|\log L(A_{e^{mt}-1}+1) - mt| > \epsilon\right), \end{aligned}$$

where we set $t := s/2$ and $m := 2n$. In the final right-hand side the first term tends to 0 because $N_t/t \to 1$ a.s., and the second does so too because of Lemma 2.9. □

**Lemma 2.12:** *As $n \to \infty$,*

$$\tilde{L}_n \circ \Phi_n = -\tilde{C}_n + o_P(1). \tag{14}$$

**Proof:** Observe (recalling that $[k] = k+1$ for integer $k$) that

$$\begin{aligned} \tilde{L}_n \circ \Phi_n &= \Big(n^{-1/2}\big(\log L(A_{N(e^{nt}-1)}+1) - A_{N(e^{nt}-1)}\big)\Big)_{0 \le t \le 1} \\ &= \Big(n^{-1/2}\big(\log L(A_{N(e^{nt}-1)}+1) - nt\big)\Big)_{0 \le t \le 1} + n^{1/2}(t - \Phi_n(t))_{0 \le t \le 1}. \end{aligned}$$

By Lemma 2.11 the first term on the right is $o_P(1)$ as $n \to \infty$. Since $\tilde{C}_n(t) = n^{1/2}\big(\Phi_n(t) - t\big)$, we conclude (14). □

**Proof of Theorem 1.2 :** From Lemma 2.12 we conclude that

$$(\tilde{T}_n \circ \Phi_n, \tilde{C}_n) = (\tilde{T}_n \circ \Phi_n, -\tilde{L}_n \circ \Phi_n) + o_P(1),$$

and then from Lemma 2.6 that

$$(\tilde{T}_n \circ \Phi_n, \tilde{C}_n) \xrightarrow{d} (\Psi_1(B), -B). \tag{15}$$

Now

$$\begin{aligned}\left(n^{-3/2}\Big(T_{n\Phi_n(t)} - \frac{1}{2}(nt)^2\Big)\right)_{0\le t\le 1} &= \left(n^{-3/2}\Big(T_{n\Phi_n(t)} - \frac{1}{2}\big(n\Phi_n(t)\big)^2\Big) + \frac{1}{2}n^{1/2}\big(\Phi_n(t)^2 - t^2\big)\right)_{0\le t\le 1}\\ &= \tilde{T}_n \circ \Phi_n - \left(\frac{1}{2}n^{1/2}\big(t - \Phi_n(t)\big)\big(t + \Phi_n(t)\big)\right)_{0\le t\le 1}\\ &= \tilde{T}_n \circ \Phi_n - \Big(n^{1/2}\big(t - \Phi_n(t)\big)t\Big)_{0\le t\le 1} + o_P(1),\end{aligned}$$

with the last step by (12). From (15), as addition is continuous in $D[0,1] \times D[0,1]$ at elements $(x, y)$ with no discontinuities in common ([13], §4), it follows that

$$\left(\tilde{C}_n, \left(n^{-3/2}\Big(T_{n\Phi_n(t)} - \frac{1}{2}(nt)^2\Big)\right)_{0\le t\le 1}\right) \xrightarrow{d} \Big(-B, \Psi_1(B) - \big(tB_t\big)_{0\le t\le 1}\Big). \tag{16}$$

Now from (4), writing $\tilde{t} := e^{nt} - 1$,

$$\begin{aligned}W(\tilde{t}) &= T_{A_{N(\tilde{t})}} + (\tilde{t} - S_{N_{\tilde{t}}})I_{N_{\tilde{t}}+1}\\ &= T_{n\Phi_n(t)} + (\tilde{t} - S_{N_{\tilde{t}}})I_{N_{\tilde{t}}+1},\end{aligned}$$

so that

$$\tilde{W}_n(t) = n^{-3/2}\big(T_{n\Phi_n(t)} - \tfrac{1}{2}(nt)^2\big) + n^{-3/2}(\tilde{t} - S_{N_{\tilde{t}}})I_{N_{\tilde{t}}+1}.$$

We will thus have

$$\tilde{W}_n \;=\; n^{-3/2}\Big(T_{n\Phi_n(t)} - \frac{1}{2}(nt)^2\Big)_{0\le t\le 1} \;+\; o_P(1) \tag{17}$$

provided we can show that

$$n^{-3/2} \sup_{0\le t\le e^n - 1} (t - S_{N_t})I_{N_t+1} \xrightarrow{P} 0.$$

We prove the equivalent assertion

$$(\log n)^{-3/2} \sup_{0\le t\le n} (t - S_{N_t})I_{N_t+1} \xrightarrow{P} 0. \tag{18}$$

Choose $\epsilon > 0$, then

$$P\Big((\log n)^{-3/2} \sup_{0\le t\le n} (t - S_{N_t}) I_{N_t+1} > \epsilon\Big) \le P(N_n > 2n) + P\big((\log n)^{-3/2} \max\{X_1,\dots,X_{2n+1}\} > \epsilon\big)$$
$$= P(N_n > 2n) + 1 - \big(1 - e^{-\epsilon(\log n)^{-3/2}}\big)^{2n+1}$$
$$= o(1) + (1 + o(1)) \frac{2n}{e^{\epsilon(\log n)^{3/2}}} \to 0.$$

This proves (18).

Combining (16) and (17), we have thus proved that

$$(\tilde{C}_n, \tilde{W}_n) \xrightarrow{d} (Z^C, Z^W) \quad \text{with} \quad Z_t^C = -B_t \quad \text{and} \quad Z_t^W = \Psi_1(B)_t - tB_t.$$

It remains to find the covariance function for the limit process $Z^W$. Let $Y_t = \Psi_1(B)_t = \int_0^t B_s\,ds$. For $s \le t$ we obtain

$$EY_sY_t \;=\; E\int_0^t\!\!\int_0^s B_u B_v\,du\,dv = \int_0^t\!\!\int_0^s \min\{u,v\}\,du\,dv \;=\; \frac{1}{2}s^2t - \frac{1}{6}s^3,$$

which together with the results

$$E(B_sY_t) = st - \frac{1}{2}s^2, \quad E(B_tY_s) = \frac{1}{2}s^2,$$

similarly derived, leads to

$$\operatorname{cov}(Z_s^W, Z_t^W) \;=\; EY_sY_t - sEB_sY_t - tEB_tY_s + stEB_sB_t \;=\; \frac{1}{3}s^3.$$

□

## 3. Complements and remarks

### 3.1.

A weak version of Theorem 1.1, with convergence in probability rather than almost sure convergence, has already been obtained in [12].

### 3.2.

For a lifetime distribution with bounded support and an atom at the supremum of its support the number of record renewals is obviously finite with probability 1. For discrete uniform distributions on $\{1, 2, \dots, d\}$ asymptotic normality was obtained in [10] for the total time spent in records as $d \to \infty$. That paper also contains some remarks on the interface of renewal theory and the theory of records, manifesting itself there in the 'perpetuity' $\sum_{n=1}^{\infty} \prod_{k=1}^{n} U_k$ with $(U_k)_{k\in\mathbb{N}}$ a sequence of independent random variables, each uniformly distributed on $(0, 1)$. A perpetuity appears in the present paper also, at the end of the proof of Lemma 2.3 and in the Lemma following, in the shape of the sequence

$$Y_n = e^{-R_n} + \sum_{i=1}^{n-1} V_i e^{-(R_n - R_i)}.$$

This propagates itself by iteration of the relation $Y_{n+1} = (Y_n + V_n)U_{n+1}$, where $Y_n$, $V_n$ and $U_{n+1}$ are mutually independent, $V_n$ has a unit exponential distribution and $U_{n+1} = e^{-\Delta_{n+1}}$ is uniform on $(0,1)$. The stationary distribution is that of $Y = \sum_{n=1}^{\infty} V_n \prod_{k=1}^{n} U_k$. Perpetuities of this type are discussed in [8]. The results of Lemma 2.3 can be seen as remarkable in that they show that the record value sequence takes most of the randomness out of the record times sequence; the randomness that 'remains' is encapsulated in the perpetuity.

### 3.3.

We could divide the time line into record and non-record pieces, as we have, not just for Poisson processes but for completely general renewal processes. How does the tail behaviour of the lifetime distribution affect the relative size of the record pieces? Above we have worked out the answer for exponential lifetime distributions, but we expect it to be possible to find analogues of Theorem 1.1 for other distributions. The 'exponential quantile function' (i.e. the monotone function that turns a standard exponential into a random variable with prescribed distribution) is likely to play a key role.

### 3.4.

Our results can also be seen in the larger context of *ordinary sampling vs. renewal sampling*, and its consequences for the theory of records. In order to explain this we start with a one-dimensional distribution function $F$ which, for simplicity, we assume to be continuous.

By an ordinary sample of size $n$ we mean a family $X_1, \ldots, X_n$ of independent random variables with distribution function $F$. Like many of the core results of probability theory the classical theory of records refers to this setup. In applications, however, the random variables $X_1, X_2, \ldots$, while still being independent and identically distributed, often arrive as the times between successive events in some counting process, so that instead of $X_1, \ldots, X_n$ we observe the associated counting process $N = (N_t)_{t\geq 0}$ over a fixed finite time interval. Obviously, this situation, which is what we mean by renewal sampling, changes the stochastic structure of the data. It is trivially true, for example, that observing $N$ over a fixed time interval $[0,t]$ has a truncation effect as the values of $F$ then matter only on this interval. Also, and more interestingly, the celebrated inspection paradox can be seen in this context: it essentially means that the distribution $\mathcal{L}(X_{N_t+1})$ of the inter-arrival time straddling $t$ differs from the original distribution of the $X$-variables.

How does renewal sampling change the stochastic structure of records? One of the felicitous aspects of the ordinary sampling case is the fact that for many results the distribution of the $X$-variables does not matter (as long as $F$ is continuous). That may no longer be true if we replace ordinary sampling by renewal sampling. Indeed, with

$$C_n := \#\{1 \leq i \leq n : X_i \geq X_n\}$$

the rank of the last variable in the sample, it is known that $C_n$ is uniformly distributed on the set $\{1, \ldots, n\}$, so that the distribution of $n^{-1}C_n$ converges to the uniform distribution on the unit interval as $n \to \infty$, whereas it was shown in [9] that this is no longer true if we consider $(N_t + 1)^{-1}C_{N_t+1}$ instead. For example, in the case of a unit-rate Poisson process, i.e. with $\mathcal{L}(X_1)$ exponential with mean 1, we still have convergence in distribution, but the limit distribution now has density $x \mapsto -\log x$ for $0 < x < 1$.

Of course, not all aspects of the record structure will change when passing from ordinary to renewal sampling. Consider for example the maximum

$$M_n := \max\{X_1, \ldots, X_n\}$$

of the first $n$ values for samples from the exponential distribution with mean 1. It is well known that $M_n - \log n$ converges in distribution as $n \to \infty$, and that the limit distribution is the Gumbel distribution, and it is straightforward to prove that the same holds for $M_{N_t} - \log N_t$ as $t \to \infty$. Nevertheless even this example can serve as the basis for a seemingly paradoxical consequence of random sub-sampling: the current maximum of the completed lifetimes is equal to the current record value, so that $M_{N_t} = R_{A_{N_t}}$ in the notation used in the previous sections. Now we know that $(R_n)_{n\in\mathbb{N}}$ is the sequence of partial sums associated with an i.i.d. sequence of standard exponentials, which implies that $(R_n - n)/\sqrt{n}$ converges in distribution to a standard normal as $n \to \infty$. Thus, sampling along $A_{N_t}$ removes the need to scale the centred variables by a factor that tends to infinity, as for the maximum a simple shift is enough to obtain a nontrivial limit distribution (of course, the point is that the shift obtained by inserting $N_t$ for $n$ is now random).

Finally, we wish to point out that the behaviour of records under renewal sampling has applications in the analysis of algorithms. In [11] a discrete variant of renewal sampling, with geometric inter-arrival times, appeared in connection with the analysis of an iterative procedure, von Neumann addition. There the length of the longest time interval between successive events essentially determines the number of iterations needed, which in turn is the main influence for the running time of the algorithm.